\documentclass[fleqn]{mat01}
\usepackage{times,mathtimy,amssymb,latexsym}
\begin{document}

\setcounter{page}{177} \firstpage{177}

\newtheorem{theore}{Theorem}
\renewcommand\thetheore{\arabic{theore}}
\newtheorem{theor}[theore]{\bf Theorem}

\newtheorem{lem}{Lemma}
\newtheorem{propo}{\rm PROPOSITION}
\newtheorem{coro}{\rm COROLLARY}

\title{Ergodic theory of amenable semigroup actions}

\markboth{Ali Ghaffari}{Ergodic theory of amenable semigroup
actions}

\author{ALI GHAFFARI}

\address{Department of Mathematics, Semnan University, Semnan,
Iran\\
\noindent E-mail: ghaffari1380@yahoo.com}

\volume{117}

\mon{May}

\parts{2}

\pubyear{2007}

\Date{MS received 7 September 2005; revised 23 January 2007}

\begin{abstract}
In this paper, among other things, we state and prove the mean
ergodic theorem for amenable semigroup algebras.
\end{abstract}


\maketitle

\section{Introduction}\label{S:intro}

Actions of amenable semigroups in connection with fixed point
conditions were considered in \cite{8}. The classical mean ergodic
and pointwise erogdic theorems are concerned with the convergence
of the sequence of Cesaro means \cite{14}. The mean ergodic
theorem for amenable locally compact groups was proved by
Greenleaf \cite{9,10} (for more on mean ergodic theorem on locally
compact group, the reader is referred to \cite{10} and \cite{14}).
For the pointwise ergodic theorem, refer \cite{14}.

Our primary concern in this paper is to prove a theorem for an
amenable locally compact semigroup $S$. Other aspects of the
relationship between amenability and ergodic theory are mentioned
in Corollaries~1, 2. A detailed historical account on the notion
of ergodicity is given in \cite{13}.

Let $S$ be a locally compact Hausdorff semitopological semigroup.
Let $M(S)$ be the Banach algebra of all bounded, regular Borel
measures on $S$ with total variation norm and convolution as
multiplication, and let $M_0(S)$ be the semigroup of all
probability measures in $M(S)$. Let $M(S)^*$ be the continuous
dual of $M(S)$, and 1 the linear functional in $M(S)^*$ such that
$\langle 1,\mu\rangle=\mu(S)$ for all $\mu$ in $M(S)$.

Recall that on $M(S)^{**}$ we define the first Arens product by
\begin{equation*}
\langle f\mu,\nu\rangle =\langle f,\mu*\nu\rangle,\quad  \langle
Ff,\mu\rangle=\langle F,f\mu\rangle,\quad \langle
GF,f\rangle=\langle G,Ff\rangle,
\end{equation*}
where $\mu,\nu\in M(S)$, $f\in M(S)^*$ and $F,G \in M(S)^{**}$.
A~linear functional $M\in M(S)^{**}$ is called a $\emph {mean}$ if
$\langle M,f\rangle\ge0$ whenever $f\ge0$ and if $\langle
M,1\rangle=1$. A~mean $M$ is \emph{topological left $($right$)$
invariant} if $\langle M,f\mu\rangle=\langle M,f\rangle$ $(\langle
M,\mu f\rangle=\langle M,f\rangle)$ for any $\mu\in M_0(S)$ and
$f\in M(S)^*$. We know that topologically left invariant means on
$M(S)^*$ have been studied by Riazi and Wong \cite{16} and by Wong
\cite{18,19,20}. We shall follow \cite{8} and \cite{16} for
definitions and terminologies not explained here.

\section{Main results}\label{S:P*}

Throughout the paper, $S$ is a locally compact Hausdorff
semitopological semigroup. A~semigroup $S$ is said to \emph{act}
on a set $X$ (from the right) if there exists a mapping
\begin{equation*}
(\xi,s)\mapsto\xi s,\ X\times S\to X
\end{equation*}
such that $\xi (st)=(\xi s)t$ whenever $s,t\in S$ and $\xi\in X$;
$S$ is a \emph{transformation semigroup} on $X$. We say that the
complex Banach space $X$ with a mapping $(\xi,s)\mapsto\xi s$
constitutes a \emph{right Banach $S$-module} if the mapping has
the following properties:\vspace{-.2pc}
\begin{enumerate}
\renewcommand\labelenumi{(\alph{enumi})}
\leftskip .1pc
\item for each fixed $s\in S$, the mapping $\xi\mapsto\xi s$ is
linear on $X$;

\item there exists a positive constant $k$ such that $\|\xi s\|\leq
k\|\xi\|$ $(\xi\in X,\ s\in S)$ and, for every $\xi\in X$, the
mapping $s\mapsto\xi s$ is continuous.\vspace{-.2pc}
\end{enumerate}

We define similarly a \emph{left dual $S$-module structure} on
$X^*$ by putting $\langle sf,\xi\rangle=\langle f,\xi s\rangle$
whenever $s\in S$, $\xi\in X$ and $f\in X^*$. For every $f\in
X^*$, $s\mapsto sf$ is a continuous mapping of $S$ into $X^*$ when
$X^*$ is equipped with the weak$^*$-topology.

Let $X$ be a right Banach $S$-module. We may give $X$ a right
Banach $M(S)$-module structure via a vector-valued integral. We
put
\begin{equation*}
\xi \mu=\int \xi s {\rm d}\mu(s),
\end{equation*}
for $\mu\in M(S)$, $\xi\in X$. Notice also that if $s\in S$ and
$\xi\in X$, then $\xi\delta_s=\xi s$. Finally, for every $\mu\in
M(S)$, $\xi\in X$ and $f\in X^*$, we define $\langle
\mu,f\xi\rangle=\langle f,\xi\mu\rangle.$

\begin{propo}\label{T:P*}$\left.\right.$\vspace{.5pc}

\noindent Let $\xi\in X$ be such that
$C_\xi=\overline{co}\{\xi\delta_s;  s\in S\}$ is weakly compact.
Then $C_\xi=\overline{\{\xi\mu; \mu\in M_0(S)\}}$.
\end{propo}

\begin{proof} It is clear that $C_\xi\subseteq \overline{\{\xi\mu;\
\mu\in M_0(S)\}}$. We will see that this inclusion is actually an
equality. To reach a contradiction, assume that some
$\xi\mu\in\{\xi\mu; \mu\in M_0(S)\}$ is not in $C_\xi$.
Theorem~3.4(b) of \cite{17} shows that there exist $f\in X^*$ and
$\gamma\in \Bbb R$ such that, for every $s\in
S$,
\begin{equation*}
\textrm{Re}\langle f,\xi s\rangle<\gamma<\textrm{Re}\langle
f,\xi\mu\rangle.
\end{equation*}
Therefore
\begin{equation*}
\textrm{Re}\langle f,\xi\mu\rangle=\int \textrm{Re}\langle f,\xi
s\rangle {\rm d}\mu(s)<\textrm{Re}\langle f,\xi\mu\rangle,
\end{equation*}
which is contradiction.\hfill $\Box$
\end{proof}

Let $\Omega$ be a subset of $M(S)$. We say that a net
$(\mu_\alpha)$ in $M_0(S)$ is \emph{asymptotically
$\Omega$-invariant} if $\|\mu*\mu_\alpha-\mu_\alpha\|\to 0$
whenever $\mu\in\Omega$. Asymptotical invariance property was
introduced by Day \cite{4}. Wong \cite{18} demonstrated that
$M(S)^*$ has a topologically left invariant mean if and only if
there exists a net $(\mu_\alpha)$ in $M_0(S)$ such that
$\|\mu*\mu_\alpha-\mu_\alpha\|\to 0$ whenever $\mu\in M_0(S)$. If
$S$ is a semi-foundation semigroup, the following conditions
\cite{7} are equivalent:
\begin{enumerate}
\renewcommand\labelenumi{(\alph{enumi})}
\leftskip .1pc
\item $M(S)^*$ has a topologically left invariant mean;

\item there is a net $(\mu_\alpha)$ in $M_0(S)$ such that, for
every compact subset $K$ of $S$, $\|\mu*\mu_\alpha-\mu_\alpha\|\to
0$ uniformly over all $\mu$ in $M_0(S)$ which are supported on
$K$.\vspace{-.2pc}
\end{enumerate}

Let $M(S)^*$ have a topologically left invariant mean, and let
$(\mu_\alpha)$ be a net in $M_0(S)$. It is easy to see that
$(\mu_\alpha)$ is asymptotically $M_0(S)$-invariant if and only if
\begin{equation*}
\lim_\alpha\|\mu*\mu_\alpha\|=|\mu(S)|
\end{equation*}
whenever $\mu\in M(S)$.

Let $X$ be a right Banach $S$-module. Let $(\mu_\alpha)$ be an
asymptotically $M_0(S)$-invariant in $M_0(S)$. For each $\alpha$,
we define $E_\alpha[\xi]\in X$ by
\begin{equation*}
E_\alpha[\xi]=\int\xi s{\rm d}\mu_\alpha(s),\quad \xi\in X.
\end{equation*}
We now concentrate on the mean ergodic theorem for a locally
compact semigroup. The following result can be referred to as `the
mean ergodic theorem' for amenable locally compact semigroups (see
\cite{10,14,15}, for details$)$.

\begin{theor}[\!]\label{T:P*}
Let $M(S)^*$ have a topologically right invariant mean. Suppose
$(\mu_\alpha)$ is an asymptotically $M_0(S)$-invariant in
$M_0(S)$.\vspace{-.2pc}
\begin{enumerate}
\renewcommand\labelenumi{\rm (\arabic{enumi})}
\leftskip .1pc
\item If $\xi\in X$ is such that
$C_\xi=\overline{co}\{\xi\delta_s; s\in S\}$ is weakly compact{\rm
,} then $(E_\alpha[\xi])$ converges to a fixed point of $C_\xi$
that is the unique fixed point in $C_\xi,$ and is therefore
independent of the choice of $(\mu_\alpha)$.

\item Suppose that $C_\xi$ is weakly compact for each $\xi\in X$.
Let $Y=\{\xi\in X;\ \xi s=\xi \text{ for all}\  s\in S\},$ and let
$Z$ be the closed vector subspace of $X$ generated by $\{\xi-\xi
s;\ s\in S,\ \xi\in X\}$. Then $X=Y\bigoplus Z$ and $(E_\alpha)$
converges strongly to the projection of $X$ onto $Y$.

\item Let $C_\xi$ be weakly compact for any $\xi\in X$. If
$P\hbox{\rm :}\ X\to Y$ is the projection onto $Y$ associated with
the direct sum decomposition of $(2),$ then $P(\xi)=Y\bigcap
C_\xi$ for all $\xi\in X$.\vspace{-.2pc}
\end{enumerate}
\end{theor}

\begin{proof}$\left.\right.$\vspace{.5pc}

\noindent (1) Let $M$ be a topologically right invariant mean on
$M(S)^*$. If $\mathcal{M}$ is the convex set of all means on
$M(S)^*$, it is well known that $\mathcal{M}$ is weak$^*$ compact
and that $M_0(S)$ is weak$^*$ dense in $\mathcal{M}$ (see
\cite{3,11} for details). Therefore there is a net $(\mu_\alpha)$
in $M_0(S)$ such that $\mu_\alpha\to M$ in the weak$^*$ topology.
A~subnet $(\xi\mu_\beta)$ of $(\xi\mu_\alpha)$ converges in the
weak topology to an element $\zeta$ in the compact subset
$C_\xi=\overline{\{\xi\mu;\mu\in M_0(S)\}}$.  For every $s\in S$
and $f\in X^*$, we have
\begin{align*}
\langle f,\zeta s\rangle &= \langle
sf,\zeta\rangle=\lim_\beta\langle sf,\xi\mu_\beta\rangle=
\lim_\beta\langle \mu_\beta(sf),\xi\rangle=
\lim_\beta\langle(\mu_\beta*\delta_s)f,\xi\rangle\\[.3pc]
&= \lim_\beta\langle\mu_\beta(\delta_sf),
\xi\rangle=\lim_\beta\langle\mu_\beta,(\delta_sf)
\xi\rangle=\lim_\beta\langle\mu_\beta,\delta_s(f\xi)\rangle\\[.3pc]
&= \langle M,\delta_s(f\xi)\rangle=\langle M,f\xi\rangle
=\lim_\beta\langle\mu_\beta,f\xi\rangle\\[.3pc]
&= \lim_\beta\langle\mu_\beta f,\xi\rangle=\lim_\beta\langle
f,\xi\mu_\beta\rangle=\langle f,\zeta\rangle.
\end{align*}
So $\zeta s=\zeta$ for every $s\in S$, that is, $\zeta$ is a fixed
point under the action of $S$. If $\epsilon>0$, we determine
$\alpha_1,\dots,\alpha_m$ and $s_1,\dots,s_m\in S$ such that
$\sum_{i=1}^m\alpha_i=1$ and $\|\xi^\prime-\zeta\|<\epsilon$,
where $\xi^\prime=\sum_{i=1}^m\alpha_i\xi s_i$. We have
\begin{align*}
\|E_\alpha[\xi]-\zeta\|&= \|\xi\mu_\alpha-\zeta\|\leq\|
\xi\mu_\alpha-\xi^\prime\mu_\alpha\|+
\|\xi^\prime\mu_\alpha-\zeta\|\\[.3pc]
&= \|\xi\sum_{i=1}^m(\alpha_i\
\mu_\alpha-\alpha_i\delta_{s_i}*\mu_\alpha)\|+
\|\xi^\prime\mu_\alpha-\zeta\mu_\alpha\|\\[.3pc]
&\leq
k\|\xi\|\sum_{i=1}^m\alpha_i\|\mu_\alpha-\delta_{s_i}*\mu_\alpha\|+k\|\xi^
\prime-\zeta\|\\[.3pc]
&<
k\|\xi\|\sum_{i=1}^m\|\mu_\alpha-\delta_{s_i}*\mu_\alpha\|+k\epsilon.
\end{align*}
Hence
\begin{equation*}
\lim_\alpha\|E_\alpha[\xi]-\zeta\|=0.
\end{equation*}
For every $\mu\in M_0(S)$, we have
$E_\alpha[\xi\mu]=\xi\mu*\mu_\alpha$. Hence
\begin{align*}
\lim_\alpha\|E_\alpha[\xi\mu]-\zeta\|&= \lim_\alpha\|
\xi\mu*\mu_\alpha-\xi\mu_\alpha\|\\[.3pc]
&\leq k\lim_\alpha\|\xi\|\|\mu*\mu_\alpha-\mu_\alpha\|=0
\end{align*}
and $\lim_\alpha\|E_\alpha[\varsigma]-\zeta\|=0$ for every
$\varsigma\in C_\xi$. Now if there exists another fixed point
$\zeta^\prime$ in $C_\xi$, then for every $\alpha$ we would have
\begin{equation*}
E_\alpha[\zeta^\prime]=\int\zeta^\prime s{\rm
d}\mu_\alpha(s)=\zeta^\prime
\end{equation*}
and
\begin{equation*}
\zeta^\prime=\lim_\alpha E_\alpha[\zeta^\prime]=\zeta.
\end{equation*}
This shows that $(E_\alpha[\xi])$ converges to the unique fixed
point $\zeta$ in $C_\xi$. Now if there exists another
asymptotically $M_0(S)$-invariant $(\nu_\beta)$ in $M_0(S)$, it is
easy to see that $\lim_\beta E_\beta[\xi]=\zeta$.

\noindent (2) Define an operator $P\hbox{\rm :}\ X\to X$ by
$P(\xi)=\lim_\alpha E_\alpha[\xi]$. Since $X$ is a right
$S$-module, it is easy to see that $\|P\|\leq k.$ For every
$\xi\in X$, $P(\xi)\in Y$ and $P(P(\xi))=\lim_\alpha
P(\xi)\mu_\alpha=P(\xi),$ $P^2=P$. If $\xi\in Y$, then
$P(\xi)=\lim_\alpha E_\alpha[\xi]=\xi\in Y$. We conclude that $Y$
is the closed subspace $P(X)$ of $X$. Let $\xi\in X$ and $s\in S$.
We have
\begin{align*}
\|P(\xi-\xi s)\| &=\|\lim_\alpha E_\alpha[\xi-\xi
s]\|=\lim_\alpha\|\xi\mu_\alpha-\xi\delta_s*\mu_\alpha\|\\[.3pc]
&\leq k\|\xi\|\lim_\alpha\|\delta_s*\mu_\alpha-\mu_\alpha\|=0.
\end{align*}
So $\xi-\xi\delta_s\in Y^\perp$. Conversely if $\zeta\in Y^\perp$,
$P(\zeta)=0$ is the unique fixed point in $C_\zeta$ and, for every
$\epsilon>0$, there exist $\alpha_1,\dots,\alpha_m$ and
$s_1,\dots,s_m$ in $S$ such that $\sum_{i=1}^m\alpha_i=1$ and
$\|\sum_{i=1}^m\alpha_i\zeta\delta_{s_i}\|<\epsilon.$ Hence
$\|\zeta- \sum_{i=1}^m\alpha_i(\zeta-\zeta\delta_{s_i})\|
=\|\sum_{i=1}^m\alpha_i\zeta\delta_{s_i}\|<\epsilon,$ such that
$\zeta\in Z$.

\noindent (3) This is trivial.\hfill $\Box$
\end{proof}

\begin{coro}\label{T:P*}$\left.\right.$\vspace{.5pc}

\noindent Let $M(S)^*$ have a topologically right invariant mean.
Suppose further that $(\mu_\alpha)$ is an asymptotically
$M_0(S)$-invariant in $M_0(S)$. Suppose that $S$ acts from the
right on a reflexive Banach space $X$. Then $(E_\alpha[\xi])$
converges to an element $\zeta$ in $X$ that is $S$-invariant.
\end{coro}

\begin{proof} For every $\xi\in X$, the set
$C_\xi=\overline{co}\{\xi s; s\in S\}$ is a weak$^*$-closed subset
of the closed ball of radius $k\|\xi\|$ in $(X^*)^*$. This ball is
weak$^*$ compact by Alaoglu's theorem and hence $C_\xi$ is weak
compact in the reflexive space $X$ (p.~111 of \cite{17}). The
statement is an immediate consequence of Theorem~1.\hfill $\Box$
\end{proof}

We recall that, for a semigroup $S$, $M_a(S)$ is the set of all
measures $\mu\in M(S)$ such that both mappings $s\mapsto
|\mu|*\delta_s$ and $s\mapsto \delta_s*|\mu|$ from $S$ into $M(S)$
are weakly continuous (see \cite{1} and \cite{5} for
definition$)$. A~semigroup $S$ is called a \emph{foundation
semigroup} if $\bigcup\{\textrm{supp}\,\mu;\ \mu\in M_a(S)\}$ is
dense in $S$. It is known that $M_a(S)$ admits a positive
approximate identity with norm 1 \cite{5}.

A functional $f$ in $M_a(S)^*$ is said to be \emph{weakly almost
periodic} if the set $\{f\mu;\ \mu\in M_a(S),\ \|\mu\|\leq 1\}$ is
relatively weakly compact. We denote by $\textrm{wap}(M_a(S))$ the
closed subspace of $M_a(S)^*$ consisting of all the weakly almost
periodic functionals in $M_a(S)^*$ (for more on weakly almost
periodic functionals, the reader is referred to \cite{2} and
\cite{12}).

\begin{coro}\label{T:P*}$\left.\right.$\vspace{.5pc}

\noindent Let $S$ be a foundation{\rm ,} locally compact Hausdorff
topological semigroup with identity. Let $M(S)^*$ have a
topologically right invariant mean. Suppose $(\mu_\alpha)$ is an
asymptotically $M_0(S)$-invariant in $M_0(S)$.
\begin{enumerate}
\renewcommand\labelenumi{\rm (\arabic{enumi})}
\leftskip .1pc
\item If $f\in{\rm wap}(M_a(S)),$ then $E_\alpha[f]$
converges to an element $f^\prime$ in ${\rm wap}(M_a(S))$ that is
$S$-invariant{\rm ,} that is{\rm ,} $f^\prime s=s$ for all $s\in
S;$ $f^\prime$ is independent of the choice of $(\mu_\alpha)$.

\item Let $M$ be the topologically left invariant mean on
${\rm wap}(M_a(S))$. Then $Y=\Bbb C1$ and
\begin{equation*}
\hskip -1.25pc Z=\{f\in{\rm wap}(M_a(S));\ \langle M,f\rangle=0\}.
\end{equation*}
\end{enumerate}
\end{coro}

Note that if $(\mu_\alpha)$ is an asymptotically
$M_0(S)$-invariant in $M_0(S)$, then
\begin{equation*}
\|\mu*\mu_\alpha-\mu_\alpha\|\to 0
\end{equation*}
for all $\mu\in M_0(S)$. We choose $\nu\in M_a(S)\bigcap M_0(S)$
and let $\nu_\alpha=\mu_\alpha*\nu$, $\alpha\in I$. Clearly
$\|\mu*\nu_\alpha-\nu_\alpha\|\to 0$ for all $\mu\in M_0(S)$. If
$M$ is a weak$^*$ cluster point of $(\nu_\alpha)$, then $M$ is a
topologically left invariant mean on $\textrm{wap}(M_a(S))$.

\begin{proof}$\left.\right.$\vspace{.5pc}

\noindent (1) We apply Theorem~1 to the Banach space $X= {\rm
wap}(M_a(S))$. If $f\in{\rm wap}(M_a(S))$ and $s\in S$, let
\begin{equation*}
\langle fs,\mu\rangle=\langle f,\delta_s*\mu\rangle,\quad \mu \in
M_a(S).
\end{equation*}
We check the continuity of the mapping $s\to fs$ for any $f\in
{\rm wap}(M_a(S))$. Let $f\in{\rm wap}(M_a(S))$ and $(e_\alpha)$
be a bounded approximate identity in $M_a(S)$. It is easy to see
that $fe_\alpha\to f$ in the weak$^*$ topology. By the compactness
of $\{f\mu;\ \mu\in M_a(S),\ \|\mu\|\leq 1\}$, we can suppose that
$fe_\alpha\to f$ in the weak topology, passing to a subnet, if
necessary. Clearly $f\in M_a(S)^*M_a(S)$, since $M_a(S)^*M_a(S)$
is a Banach space and $fe_\alpha\in M_a(S)^*M_a(S)$ for all
$\alpha\in I$. By Lemma~2.2 in \cite{6}, $s\mapsto fs$ is
continuous. Obviously $X$ is a right Banach $S$-module.

If $f\in{\rm wap}(M_a(S))$, then $\{fs;\ s\in S\}$ is relatively
weakly compact. So Theorem~1 asserts that, for every $f\in X$,
$(E_\alpha[f])$ converges to an element $f^\prime$ in $X$ that is
fixed; $f^\prime s=f^\prime$ and $\langle
f^\prime,\delta_s*\mu\rangle=\langle f^\prime,\mu\rangle$ for
every $s\in S$ and every $\mu\in M_a(S)$.

\noindent $(2)$ By Theorem~1, $Y$ is the set of all $f\in X$ such
that $fs=f$ for all $s\in S$. If $f\in Y$, then $\langle
f,\delta_s*\mu\rangle=\langle f,\mu\rangle$ for all $s\in S$ and
$\mu\in M_a(S)$. Consequently, $f=\langle f,\delta_e\rangle 1$.

We know that $Z$ admits a dense subset consisting of functionals
of the form $g=\Sigma_1^n f_i-f_is_i$, where $f_i\in{\rm
wap}(M_a(S))$ and $s_i\in S$. By the definition of $M$, we have
$\langle M,g\rangle=0$ and so $Z\subseteq\{f\in{\rm wap}(M_a(S));\
\langle M,f\rangle=0\}$. For the reverse inclusion, let $f\in{\rm
wap}(M_a(S))$ and $\langle M,f\rangle=0$. By Theorem~1, we can
write $f=f_1+f_2$ for some $f_1\in Y$ and $f_2\in Z$, and so
\begin{align*}
0&= \langle M,f\rangle=\langle M,f_1\rangle+\langle
M,f_2\rangle\\[.3pc]
&= \langle M,f_1\rangle+0=\langle M,f_1\rangle,
\end{align*}
since $Z\subseteq\{f\in \textrm{wap}(M_a(S));\ \langle
M,f\rangle=0\}.$ Consequently $\langle M,f_1\rangle=0$. For each
$\mu\in M_0(S)$, we have
\begin{align*}
0&= \langle M,f_1\rangle=\langle
M,f_1\mu\rangle=\lim_\alpha\langle\mu_\alpha,f_1\mu\rangle\\[.3pc]
&= \lim_\alpha \langle
f_1,\mu*\mu_\alpha\rangle=\lim_\alpha\langle
f_1\mu*\mu_\alpha,\delta_e\rangle\\[.3pc]
&= \lim_\alpha\langle f_1\mu,\delta_e\rangle=\langle
f_1,\mu\rangle.
\end{align*}
Therefore $f_1=0$. It follows that $Z=\{f\in\textrm{wap}(M_a(S));\
\langle M,f\rangle=0\}$.\hfill $\Box$
\end{proof}

\section*{Acknowledgment}

I would like to thank the referee for his/her careful reading of
my paper and many valuable suggestions.

\end{document}